\documentclass[10pt,twoside]{siamltex}
\usepackage{amsfonts,epsfig}
\usepackage{hyperref}

\usepackage{graphics}
\usepackage{epsfig}
\usepackage{amssymb,amsmath,mathrsfs}
\usepackage{amsfonts}
\usepackage{color}
\usepackage{mathtools}
\usepackage{algorithmic}
\usepackage{algorithm}
\usepackage{textcomp}

\usepackage{tikz}
\usetikzlibrary{arrows}
\usetikzlibrary{decorations.markings}
\tikzstyle{block}=[draw opacity=0.7,line width=1.4cm]
\tikzset{
  big black arrow/.style={
    decoration={markings,mark=at position 1 with {\arrow[scale=2.5,black]{>}}},
    postaction={decorate},
    shorten >=0.4pt},
    line/.style={draw, ->}}

\numberwithin{equation}{section}

\newcommand {\rref}{\mathop{\rm rref}\nolimits}
\newcommand {\tril}{\mathop{\rm tril}\nolimits}

\begin{document}

%  Leave these commented lines here
% \input{elaheader-volx-xx.tex}
% \setcounter{page}{1}

% \renewcommand{\thefootnote}{\fnsymbol{footnote}}
% \renewcommand{\thefootnote}{\arabic{footnote}}
% \renewcommand{\theequation}{\thesection.\arabic{equation}}

\bibliographystyle{plain}
\title{
Generalized Jacobi and Gauss-Seidel method for solving non-square linear systems%
%\thanks{Received by the editors on Month x, 200x.Accepted for publication on Month y, 200y   Handling Editor: .}
}
% Leave blank; editors will write the exact dates above

\author{
Manideepa Saha\thanks{Department of Mathematics, National Institute of Technology Meghalaya, Shillong 793003, India
(manideepa.saha@nitm.ac.com).}
% Remember to put \and between any two authors
%\and Saikat Mukherjee\thanks{Department of Mathematics, National Institute of Technology Meghalaya, Shillong 793003, India (saikat.mukherjee@nitm.ac.in).}
}
% Note that \footnotemark[3]} is used for the third author
% because of the same affiliation for the second and third authors.
% If the same affiliation is to be used for the first and second authors,
% \footnotemark[2] should be used instead of \thanks{} for the second author.

% Authors and running title to go on top of each page
\pagestyle{myheadings}
\markboth{M.\ Saha
%and S.Mukherjee
 }{Generalized Jacobi and Gauss-Seidel method }
\maketitle

\begin{abstract} The main goal of this paper is to generalize Jacobi and Gauss-Seidel methods for solving non-square linear system. Towards this goal, we present iterative procedures to obtain an approximate solution for  non-square linear system. We derive sufficient conditions for the convergence of such iterative methods. Procedure is given to show that how an exact solution can be obtained from these methods. Lastly, an example is considered to compare these methods with other available method(s) for the same.

%are among the most stationary iterative methods for solving square linear system. In this paper, we describe a generalization of both Jacobi and Gauss-Seidel methods to non-square linear system. .

\end{abstract}


%%%%%%%%%%%%%%%%%%%%%%%%%%%%%%%%%%%%%%%%%%%%%%%%%%%%%%%%%%%%%

\section{Introduction} \label{intro3} 
The idea of solving large square systems of linear equations by iterative methods is certainly not new, dating back at least to Gauss [1823]. Jacobi and Gauss-Seidel methods are most stationary iterative methods that date to the late eighteenth century, but they find current application in problems where the matrix is sparse. Both these methods are designed for finding an approximate solution to  square linear systems. The purpose of this paper is to provide an iterative procedure for solving a non-square linear system. In particular, we generalize Jacobi and Gauss-Seidel methods for non-square linear system.

Very recently, in \cite{WheA17}, authors considered non-square linear system with number of variables is more than that of equations and described a new iterative procedure along with a convergence analysis. They also discussed some of its applications. Motivated by their work, and using similar technique, we generalize Jacobi and Gauss-Seidel procedure for solving non-square linear system with number of variables is more than that of equations. More specifically, we apply Jacobi or Gauss-Seidel iterative method for the square part of the system and apply the iterative method described in \cite{WheA17} for the non-square part of the system to obtain an approximate solution of the system. We also derive sufficient conditions for the convergence of such procedure. Finally, we provide a procedure to obtain an exact solution of the system . A numerical illustration has been given for the same, and to compare the procedure with the procedure available in~\cite{WheA17}.

\section{Preliminaries} \label{convention} Consider the linear system 
$Ax=b$, where $A\in\mathbb{R}^{m,n}$, and $b\in\mathbb{R}^{m}$ with $m<n$. We now state iterative procedure described in \cite{WheA17}. If $z^{(k)}$ is the approximate solution obtained after $k$th iteration, then 
\begin{equation}\label{eqnWA1}
z^{(k+1)}=z^{(k)}+s(A)d^{k}_{A}
\end{equation}
 where $s(A)$ is the $n\times m$ matrix with $1,0,-1$ as its elements, and they indicate the signs of elements of $A^{T}$, and $d^{(k)}_{A}$ is the column $m$-vector of weighted residuals, that is, $d^{k}_{A}=\left[\frac{b_{1}-A_{1}z^{(k)}}{m\|A_{1}\|_{1}},\frac{b_{2}-A_{2}z^{(k)}}{m\|A_{2}\|_{1}},\ldots,\frac{b_{1}-A_{m}z^{(k)}}{m\|A_{m}\|_{1}}\right]^{T}$, where $A_{i}$ denotes the $i$th row vector of $A$ and $\|.\|_{1}$ denotes the $l_{1}$-norm.

%

%%%%%%%%%%%%%%%%%%%%%%%%%%%%%%%%%%
\section{Generalized Jacobi and Gauss-Seidel Method} In this section we describe iterative procedures, based on the iterative procedure (\ref{eqnWA1}) and on Jacobi or Gauss-Seidel procedure for square matrix, and called them as generalized Jacobi or generalized Gauss-Seidel method, respectively.

\subsection{Steps for generalized Jacobi and generalized Gauss-Seidel Method}
 Our notation will be similar to that of used in the previous section. As $m>n$, we can partition $A$ as $A=[B~ \tilde{B}]$ with $B\in\mathbb{R}^{m,m}$, and $\tilde{B}\in\mathbb{R}^{m,n-m}$. Analogously, partition the initial approximation $x^{(0)}=\left[\begin{array}{l}
  x^{(0)}_{1}\\x^{(0)}_{2}
  \end{array}\right]$, conformably with $A$. Let $\epsilon$ be a given threshold. Followings are the steps of the procedures.
  
  If $\|Ax^{(0)}-b\|<\epsilon$, then initial guess $x^{(0)}$ will be an approximate solution of the system $Ax=b$. Otherwise, following steps will be executed:
  
 \begin{description}
 \item[Step 1 :] Set iteration count $k\leftarrow 0$. Calculate $\tilde{b}^{(k)}=b-Bx^{(k)}_{1}$. Apply iteration procedure proposed in~\cite{WheA17} to the system $\tilde{B}y=\tilde{b}^{(k)}$ with initial guess $x^{(k)}_{2}$ to get $x^{(k+1)}_{2}$, that is,  
 \[x^{(k+1)}_{2}\leftarrow x^{(k)}_{2}+s(\tilde{B}) d^{(k)}\]
 
  where  $s(\tilde{B})$ is an $n-m\times m$ matrix, elements of which represents the signs of elements of $\tilde{B}^T$, and $d^{(k)}$ is an $m$-vector whose $i$th entry is $d^{(k)}_{i}=\frac{\tilde{b}^{(k)}_{i}-\tilde{B}_{i}x^{(k)}_{2}}{m\|\tilde{B}_{i}\|_{1}}$.
 
 \item[Step 2 :] Set $\hat{b}^{(k)}=b-\tilde{B}x^{(k+1)}_{2}$. Apply Jacobi or, Gauss-Seidel method to the square linear system $Bz=\hat{b}^{(k)}$ with initial guess $x^{(k)}_{1}$ to get $x^{(k+1)}_{1}$, that is,
  \[x^{(k+1)}_{1}\leftarrow D^{-1}(-Rx^{(k)}_{1}+\hat{b}^{(k)})~~(\text{ in case Jacobi method})\]
 where $D=\diag(B)$, and $R=B\setminus D$,
 \[\text{ or, }x^{(k+1)}_{1}\leftarrow L^{-1}(-Rx^{(k)}_{1}+\hat{b}^{(k)})~~(\text{ in case Gauss-Seidel method})\]
 where $L=\tril(B)$, and $R=B\setminus L$.
 
 \item[Step 3 :] Set $x^{(k+1)}=\left[\begin{array}{l}
  x^{(k+1)}_{1}\\x^{(k+1)}_{2}
  \end{array}\right]$, where $x^{(k+1)}_{1}$ is obtained from Step I and $x^{(k+1)}_{2}$ is obtained from Step II.
  
 \item[Step 4 :] Repeat the above steps untill $\|Ax^{(k+1)}-b\|<\epsilon .$
 
\end{description}

\subsection{Convergence analysis} Assume that the system $Ax=b$ has a solution and let $x^{(k)}=\left[\begin{array}{l}
  x^{(k)}_{1}\\x^{(k)}_{2}
  \end{array}\right]$ be an approximation solution obtained after $k$th iteration step, where $x^{(k)}_{1}\in\mathbb{R}^{m}$ and $x^{(k)}_{2}\in\mathbb{R}^{n-m}$.

From Step I of the iterative procedure we get, 
\begin{equation}\label{eqn0}
x^{(k+1)}_{2}= x^{(k)}_{2}+s(\tilde{B}) d^{(k)}
\end{equation} 

%\begin{equation}\label{eqn1}
%\tilde{B}x^{(k+1)}_{2}= \tilde{B}x^{(k)}_{2}+\tilde{B}*s(\tilde{B}) d^{(k)}
%\end{equation}
%where * is  the matrix multiplication.

If we use Jacobi method in step II, then  
\begin{equation}\label{eqn1}x^{(k+1)}_{1}=D^{-1}\left(-Rx^{(k)}_{1}+\hat{b}^{(k)}\right)\end{equation}

 Note that $d^{(k)}_{i}=\frac{\tilde{b}^{(k)}_{i}-\tilde{B}_{i}x^{(k)}_{2}}{m\|\tilde{B}_{i}\|_{1}}=\frac{b_{i}-B_{i}x^{(k)}_{1}-\tilde{B}_{i}x^{(k)}_{2}}{m\|\tilde{B}_{i}\|_{1}}=\frac{b_{i}-A_{i}x^{(k)}}{m\|\tilde{B}_{i}\|_{1}}$.
 If we take $N(\tilde{B})=\diag(\|\tilde{B}_{1}\|_1, \|\tilde{B}_{2}\|_1,\ldots, \|\tilde{B}_{m}\|_1)$, then $d^{(k)}=\frac{1}{m}N(\tilde{B})^{-1}(b-Ax^{(k)})$, if $N(\tilde{B})$ is nonsingular. Hence equation (\ref{eqn1}) can be further reduced to
\begin{flalign}\label{eqn2}
x^{(k+1)}_{1}=x^{(k)}_{1}-D^{-1}\left(I-\frac{1}{m}\tilde{B}*s(\tilde{B})*N(\tilde{B})^{-1}\right)\left(Ax^{(k)}-b\right) 
\end{flalign}
where * denotes the matrix multiplication. As $Ax^{(k+1)}=Bx^{(k+1)}_{1}+Bx^{(k+1)}_{2}$, equations (\ref{eqn0}) and (\ref{eqn2}) imply that
 \begin{flalign} \label{eqn3}
 Ax^{(k+1)}-b&= (I-BD^{-1})\left(I-\frac{1}{m}\tilde{B}*s(\tilde{B})*N(\tilde{B})^{-1}\right)\left(Ax^{(k)}-b\right)
 \end{flalign}
 
This shows that if  $\|I-BD^{-1}\|<1$ and $\|mI-\tilde{B}*s(\tilde{B})*N(\tilde{B})^{-1}\|<m$, then $x=\lim\limits_{k}x^{(k)}$, if exists, is a solution of the system $Ax=b$. 

 Assume that $\|I-BD^{-1}\|<1$ and $\|mI-\tilde{B}*s(\tilde{B})*N(\tilde{B})^{-1}\|<m$. We now the existance of $x$. From (\ref{eqn0}) and (\ref{eqn2}) we get
{\small\begin{flalign}\label{eqn4}
  \|x^{(k+1)}-x^{(k)}\|_1 &= \|x^{(k+1)}_{1}-x^{(k)}_{1}\|_1+\|x^{(k+1)}_{2}-x^{(k)}_{2}\|_1 \nonumber\\
  %&=\|D^{-1}\left(Ax^{(k)}-b+\tilde{B}*s(\tilde{B}) d^{(k)}\right)\|_1+\|s(\tilde{B}) d^{(k)}\|_1\nonumber\\
   &=\|D^{-1}\left(I-\frac{1}{m}\tilde{B}*s(\tilde{B})*N(\tilde{B})^{-1}\right)(b-Ax^{(k)})\|_1+\|\frac{1}{m}s(\tilde{B})*N(\tilde{B})^{-1}(b-Ax^{(k)})\|_1\nonumber\\
   &\leq\left[\|D^{-1}\left(I-\frac{1}{m}\tilde{B}*s(\tilde{B})*N(\tilde{B})^{-1}\right)\|_1+\frac{1}{m}\|s(\tilde{B})*N(\tilde{B})^{-1}\|_1\right]\|Ax^{(k)}-b\|_1
 \end{flalign}}

From equations (\ref{eqn2}) and (\ref{eqn4}), it can be easily verified that $\{x^{k}\}$ is a Cauchy-sequence, and hence  $x=\lim\limits_{k}x^{(k)}$ exists.

From the above discussion we can coclude that the our proposed generalization of Jacobi method will converge if $\|I-BD^{-1}\|<1$ and $\|mI-\tilde{B}*s(\tilde{B})*N(\tilde{B})^{-1}\|<m$, for any matrix norm $\|.\|$. 

Similarly, it can be proved that if we apply Gauss-Seidel method in Step II, the generalized  Gauss-Seidel method will converge if $\|I-BL^{-1}\|<1\|$ and $\|mI-\tilde{B}*s(\tilde{B})*N(\tilde{B})^{-1}\|<m$.

\section{Application and Comparison} Obtaining an exact solution of the system $Ax=b$. 

Let $[\bar{A}~\bar{b}]$ be the reduced row echelon form of $[A~b]$, and let $\rank(A)=m$. We now apply our procedure on the system $\bar{A}x=\bar{b}$.

As $\rank(A)=m$, $B=I$. Then $x^{(k+1)}_{2}=x^{k}_{2}+s(\tilde{B})d^{(k)}$, and $x^{(k+1)}_{1}=\hat{b}^{(k)}$, because in both procedures $R=0$. Thus $Ix^{(k+1)}_{1}=b-\tilde{B}x^{(k+1)}_{2}$, and hence $\bar{A}x^{(k+1)}=b$. Hence in each iteration we obtain an exact solution of the system.

We now compare our method with the procedure described in \cite{WheA17} with numerical illustration. Consider the system $Ax=b$, with
\[[A~b]=\left[\begin{array}{rrrrrrrrrr}
2 & 4& -3 & 1 & 0 &5 &-7 &8 &| &38\\
3 &2 &10 &-4 &-1 &-6 &4 &1 &| &20\\
9 &7 &3 &2 &0 &0 &-4 &2 &| &39\\
6 &4 &0 &-1 &-1 &3 &10 &5 &| &-16\\
5 &2 &-3 &-7 &-5 &4 &8 &-8 &|&-30
\end{array}\right]\]
so that
\[[\bar{A}~\bar{b}]=\rref([A~b])=\left[\begin{array}{rrrrrrrrrr}
1 & 0& 0 & 0 & 0 &0.2 &5.7 &0.6 &| &-20.5\\
0 &1 &0 &0 &0 &0.6 &-4.2 &2.0 &| &16.2\\
0 &0 &1 &0 &0 &0 &-4 &2 &| &9.1\\
0 &0 &0 &1 &0 &3 &10 &5 &| &41.2\\
0 &0 &0 &0 &1 &4 &8 &-8 &|&-83.1
\end{array}\right]\]
Our iterative procedure (generalized Jacobi) method applying to the system $\bar{A}x=\bar{b}$ with initial approximation $x^{0}=[2,~0,~-1,~2,~0,~0,~-3,~1]^{T}$ resulted to a non-basic approximation solution
$x=[2.9734,~2.2736,~-4.0025,~-2.9777,~-1.9033,~-1.5452,~-3.9452,~-0.7452]^{T}$.
Applying the method in \cite{WheA17} to the system $\bar{A}x=\bar{b}$ with the same initial guess $x^{0}$ gives the resutant approximate solution
$y=[2,~0,~-1,~-2,~-1.3015,~-1.3015,~-4.3015,~-0.3015]^{T}$.

 The residual obtained for the method in \cite{WheA17} is $\|b-Ay\|_{1}=11.7462$, whereas our method resulting to an exact non-basic solution $x$.

\section{Conclusion} We have described a possible generalization of Jacobi and Gauss-Seidel method to non-square linear system. We also discussed the convergence of the proposed procedures. It is shown that an exact solution of a system $Ax=b$ can be obtained if we apply our methods to its reduced row echelon system, and a comparision has been done with the method discussed in \cite{WheA17}.

We now conclude this paper by listing some ideas that we did not pursue, which however may lead to further progress.

We know that both Jacobi and Gauss-Seidel method converges for diagonally dominant coefficient matrices. One may check if the proposed methods converges for row diagonally coefficient matrices.

As SOR method (Successive overrelaxation) is also one of the most stationary iterative method for solving linear square system, one may use a similar technique to generalize SOR method, and check its convergence criteria.


\begin{thebibliography}{1}
\bibitem{WheA17}
I.~Wheaton and S.~Awoniyi.
\newblock  {\em A new iterative method for solving non-square systems of linear equations}.
\newblock  Journal of Computaional and Applied Mathematics , 322:1--6, 2017.


\bibitem{Hac94}
W. Hackbusch
\newblock  {\em Iteraive solution of large Sparse systems of equations}.
\newblock  Springer-Verlag, 1994.

\bibitem{Hac94}
R.S. Varga
\newblock  {\em Matrix iteraive analysis}.
\newblock  Springer, 2000.

\bibitem{Hac94}
M.T. Micheal
\newblock  {\em Scientific Computing: An introductory survey}.
\newblock  McGraw Hill, 2002.

\bibitem{Hac94}
C.T. Kelley
\newblock  {\em Iterative methods for linear and nonlinear equations}.
\newblock  SIAM, Philadelphia 1995.






\end{thebibliography}
\end{document}